\documentclass[a4paper,11pt]{article}
\usepackage{latexsym}
\usepackage{amssymb}
\usepackage{enumerate}
\usepackage{amsfonts}
\usepackage{amsmath}
\usepackage[dvips]{graphicx}
\usepackage{epsf}
\newcommand{\qed}{$\Box$}

\newenvironment{@abssec}[1]{%
    \if@twocolumn

      \section*{#1}%
    \else

      \vspace{.05in}\footnotesize
      \parindent .2in
 {\upshape\bfseries #1. }\ignorespaces
    \fi}

    {\if@twocolumn\else\par\vspace{.1in}\fi}

\newenvironment{keywords}{\begin{@abssec}{\keywordsname}}{\end{@abssec}}

\newenvironment{AMS}{\begin{@abssec}{\AMSname}}{\end{@abssec}}

\newcommand\keywordsname{Key words}
\newcommand\AMSname{AMS subject classifications}
\newcommand\AMname{AMS subject classification}
\newtheorem{theorem}{Theorem}
 \newtheorem{lemma}[theorem]{Lemma}
 
 \newtheorem{proposition}[theorem]{Proposition}

\def\qed{\vbox{\hrule height0.6pt\hbox{%
  \vrule height1.3ex width0.6pt\hskip0.8ex
  \vrule width0.6pt}\hrule height0.6pt
 }}

\title{Interaction between fast diffusion and geometry of domain
\thanks{This research was partially supported by Grants-in-Aid
for Scientific Research (B) ($\sharp$ 20340031 and $\sharp$ 26287020) of
Japan Society for the Promotion of Science.}}

\author{Shigeru Sakaguchi\thanks{Research Center for Pure and Applied Mathematics,
Graduate School of  Information Science, Tohoku
University, Sendai, 980-8579,  Japan.
({\tt sigersak@m.tohoku.ac.jp}).}}
\begin{document}

\maketitle

\begin{abstract}
 Let  $\Omega$ be a domain in $\mathbb R^N$, where $N \ge 2$ and $\partial\Omega$ is not necessarily bounded. We consider  two fast diffusion equations  $\partial_t u= \mbox{\rm div}(|\nabla u|^{p-2}{\nabla u})$ and $\partial_t u= \Delta u^{m}$, where $1<p<2$ and $0<m<1$. Let $u=u(x,t)$ be the solution of either the initial-boundary value problem over $\Omega$, where the initial value equals zero and the boundary value is a positive continuous function,  or 
the Cauchy problem where the initial datum equals a nonnegative continuous function multiplied by the characteristic function of the set $\mathbb R^N\setminus \Omega$.  Choose  an open ball $B$  in $\Omega$ whose closure intersects $\partial\Omega$ only at one point, and let $\alpha >  \frac {(N+1)(2-p)}{2p}$ or $\alpha >  \frac {(N+1)(1-m)}{4}$. Then, 
we derive asymptotic estimates for the integral of  $u^\alpha$ over $B$ for short times in terms of  principal curvatures of $\partial\Omega$ at the point, 
which tells us about the interaction between fast diffusion and geometry of domain.
\end{abstract}




\pagestyle{plain}
\thispagestyle{plain}
\markboth{S. SAKAGUCHI}{Interaction between fast diffusion and geometry of domain}

\pagestyle{plain}
\thispagestyle{plain}

\section{Introduction}
\label{introduction}

\noindent
Let  $\Omega$ be a domain in $\mathbb R^N$, where $N \ge 2$ and $\partial\Omega$ is not necessarily bounded. We consider  two fast diffusion equations of the forms  $\partial_t u= \mbox{\rm div}(|\nabla u|^{p-2}{\nabla u})$ and $\partial_t u= \Delta u^{m}$, where $1<p<2$ and $0<m<1$. Let $f \in C^{0}(\partial\Omega)$ be a function satisfying 
\begin{equation}
\label{f is bounded from above and below}
0 < c_{1}\le f(x) \le c_{2}\ (x \in \partial\Omega)
\end{equation}
 for two positive constants $c_{1}$ and $c_{2}$, and let $g \in C^{0}(\mathbb R^{N})$ be a function satisfying
 \begin{equation}
 \label{g is bounded from above}
 0 \le g(x) \le c_{3}\ (x \in \mathbb R^{N})
 \end{equation}
  for a positive constant $c_{3}$. Consider the bounded solution $u=u(x,t)$  of either
  the  initial-boundary value problem:
  \begin{eqnarray}
&\partial_t u=\mbox{\rm div}(|\nabla u|^{p-2}{\nabla u})\ \ &\mbox{in }\ \Omega\times (0,\infty),\label{diffusion pLaplace}\\
&u=f\ \ &\mbox{on }\ \partial\Omega\times (0,\infty),\label{dirichlet pLaplace}\\
&u=0\ \ &\mbox{on }\ \Omega\times \{0\},\label{initial pLaplace}
\end{eqnarray}
or the Cauchy  problem:
\begin{equation}
\label{cauchy pLaplace}
\partial_t u=\mbox{\rm div}(|\nabla u|^{p-2}{\nabla u})\ \mbox{ in }\ \mathbb R^N \times (0, \infty)\quad\mbox{ and }\ u = g\mathcal X_{\Omega^c}\ \mbox{ on }\ \mathbb R^N \times \{0\},
\end{equation}
where  $\mathcal X_{\Omega^c}$ is the characteristic function of the set $\Omega^c = \mathbb R^N \setminus \Omega$. The first theorem tells us about the interaction between fast diffusion and geometry of domain for  $\partial_t u= \mbox{\rm div}(|\nabla u|^{p-2}{\nabla u})$.


\begin{theorem}
\label{th:interaction curvatures pLaplace}  Let $u$ be the solution of 
either problem {\rm (\ref{diffusion pLaplace})-(\ref{initial pLaplace})} or problem {\rm (\ref{cauchy pLaplace})}.
Let  $\alpha >  \frac {(N+1)(2-p)}{2p}$ and $x_0 \in \Omega$. Assume that the open ball $B_R(x_0)$ centered at $x_0$ and with radius $R>0$ 
is contained in  $\Omega$ and such that $\overline{B_R(x_0)} \cap \partial\Omega = \{ y_0 \}$ for some $y_0 \in \partial\Omega$ and $\partial\Omega \cap B_{\delta}(y_0)$ is of class $C^{2}$ for some $\delta >0$.
Suppose that $g(y_{0}) > 0$ for problem {\rm (\ref{cauchy pLaplace})}.
Then we have:
\begin{equation}
\label{asymptotics and curvatures pLaplace}
\lim_{t\to 0^+}t^{-\frac{N+1}{2p}}\!\!\!\int\limits_{B_R(x_0)}\! (u(x,t))^{\alpha}\ dx=
c\left\{\prod\limits_{j=1}^{N-1}\left[\frac 1R - \kappa_j(y_0)\right]\right\}^{-\frac 12}.
\end{equation}
Here, 
$\kappa_1(y_0),\dots,\kappa_{N-1}(y_0)$ denote the principal curvatures of $\partial\Omega$ at $y_0$ with 
respect to the inward normal direction to $\partial\Omega$  
and $c$ is a positive constant depending only on $p, \alpha, N,$ and either $f(y_{0})$ or $g(y_{0})$.
When $\kappa_j(y_0) = \frac 1R$ for some $j \in \{ 1, \cdots, N-1\}$, 
the formula \eqref{asymptotics and curvatures pLaplace} holds by setting the right-hand side to $\infty$ {\rm (}notice that 
$\kappa_j(y_0) \le \frac1R$
for every $j \in \{ 1, \cdots, N-1\}$ {\rm )}.
\end{theorem}

Concerning $\partial_t u= \Delta u^{m}$ with $0<m<1$, let $u=u(x,t)$ be the bounded nonnegative solution of either
  the  initial-boundary value problem: 
\begin{eqnarray}
&\partial_t u=\Delta u^{m}\ \ &\mbox{in }\ \Omega\times (0,\infty),\label{diffusion mPorous}\\
&u=f\ \ &\mbox{on }\ \partial\Omega\times (0,\infty),\label{dirichlet mPorous}\\
&u=0\ \ &\mbox{on }\ \Omega\times \{0\},\label{initial mPorous}
\end{eqnarray}
or the Cauchy  problem:
\begin{equation}
\label{cauchy mPorous}
\partial_t u=\Delta u^{m}\ \mbox{ in }\ \mathbb R^N \times (0, \infty)\quad\mbox{ and }\ u = g\mathcal X_{\Omega^c}\ \mbox{ on }\ \mathbb R^N \times \{0\}.
\end{equation}
The second theorem tells us about the interaction between fast diffusion and geometry of domain for $\partial_t u= \Delta u^{m}$.


\begin{theorem}
\label{th:interaction curvatures mPorous} Let $u$ be the solution of 
either problem {\rm (\ref{diffusion mPorous})-(\ref{initial mPorous})} or problem {\rm (\ref{cauchy mPorous})}.
Let  $\alpha >  \frac {(N+1)(1-m)}{4}$ and $x_0 \in \Omega$. Assume that the open ball $B_R(x_0)$ centered at $x_0$ and with radius $R>0$ 
is contained in  $\Omega$ and such that $\overline{B_R(x_0)} \cap \partial\Omega = \{ y_0 \}$ for some $y_0 \in \partial\Omega$ and $\partial\Omega \cap B_{\delta}(y_0)$ is of class $C^{2}$ for some $\delta >0$.
Suppose that $g(y_{0}) > 0$ for problem {\rm (\ref{cauchy mPorous})}.
Then we have:
\begin{equation}
\label{asymptotics and curvatures mPorous}
\lim_{t\to 0^+}t^{-\frac{N+1}{4}}\!\!\!\int\limits_{B_R(x_0)}\! (u(x,t))^{\alpha}\ dx=
c\left\{\prod\limits_{j=1}^{N-1}\left[\frac 1R - \kappa_j(y_0)\right]\right\}^{-\frac 12}.
\end{equation}
Here, 
$\kappa_1(y_0),\dots,\kappa_{N-1}(y_0)$ denote the principal curvatures of $\partial\Omega$ at $y_0$ with 
respect to the inward normal direction to $\partial\Omega$  
and $c$ is a positive constant depending only on $m, \alpha, N,$ and either $f(y_{0})$ or $g(y_{0})$. 
When $\kappa_j(y_0) = \frac 1R$ for some $j \in \{ 1, \cdots, N-1\}$, 
the formula \eqref{asymptotics and curvatures mPorous} holds by setting the right-hand side to $\infty$. 
\end{theorem}

When $p > 2, m> 1, \alpha = 1$, and $f\equiv g \equiv 1$, the same formulas \eqref{asymptotics and curvatures pLaplace} and \eqref{asymptotics and curvatures mPorous} were obtained for problems {\rm (\ref{diffusion pLaplace})-(\ref{initial pLaplace})} and {\rm (\ref{diffusion mPorous})-(\ref{initial mPorous})} in \cite{MSprseA2007}.  With the aid of the techniques employed in \cite{MSjde2012}, one can easily see that the formulas \eqref{asymptotics and curvatures pLaplace} and \eqref{asymptotics and curvatures mPorous} also hold true for problems {\rm (\ref{cauchy pLaplace})} and {\rm (\ref{cauchy mPorous})}. Moreover, in \cite{MSjde2012}, the nonlinear diffusion equation of the form $\partial_t u = \Delta \phi(u)$ where $\delta_{1} \le \phi^\prime(s)\le \delta_{2}\ (s \in \mathbb R)$ for some positive constants $\delta_{1}$ and $ \delta_{2}$ was also dealt with. By a little more  observation, we see that any $\alpha > 0$ is OK for these cases.

 In Theorems \ref{th:interaction curvatures pLaplace} and \ref{th:interaction curvatures mPorous}, if $p$ is close to $1$ or if $N \ge 4$ and $m$ is close to $0$, then $\alpha= 1$ can not be chosen. Indeed, when $ \alpha = \frac {(N+1)(2-p)}{2p}$ or 
$ \alpha = \frac {(N+1)(1-m)}{4}$, $c = \infty$.

The main ingredients of the proofs  of the formulas \eqref{asymptotics and curvatures pLaplace} and \eqref{asymptotics and curvatures mPorous} consist of two steps. One is the reduction to the case where $\partial\Omega$ is bounded and of class $C^2$, and where both $f$ and $g$ are constant, with the aid of the comparison principle.
The other is the construction of appropriate super- and subsolutions to the problems near $\partial\Omega$ in a short time. In fact,  in \cite{MSprseA2007}, such barriers were constructed in a set $\Omega_{\rho} \times (0,\tau],$  with
\begin{equation}
\label{neighborhood of boundary}
\Omega_{\rho} =\{ x \in \Omega\ :\ \mbox{ dist}(x,\partial\Omega) < \rho \},
\end{equation}
where $\rho$ and $\tau$ were chosen sufficiently small.  When $p > 2$ or $m> 1$, the property of  finite speed of propagation of disturbances from rest yields that
both the solution $u$ and the barriers equal zero on $\Gamma_\rho\times (0,\tau],$ where 
\begin{equation}
\label{level surfaces of distance functions}
\Gamma_{\rho} =\{ x \in \Omega\ :\ \mbox{ dist}(x,\partial\Omega) =\rho \}.
\end{equation}
This property does not occur when $1 < p < 2$ or $0 < m <1$, because of the property of  {\it infinite} speed of propagation of disturbances from rest.  Also in \cite{MSjde2012},   the equation $\partial_t u = \Delta \phi(u)$ has the property of  infinite speed of propagation of disturbances from rest.  To compare the solution with the barriers on $\Gamma_{\rho}\times (0,\tau]$, in  \cite{MSjde2012},  the result of 
Atkinson and Peletier \cite{AParma1974} concerning the asymptotic behavior of 
one-dimensional similarity solutions and the following short time behavior of $u$ obtained by \cite{MSpoincare2010} play a key role:
\begin{equation}
\label{varadhan formula}
\lim_{t \to 0^+} -4t\Phi(u) = \mbox{ dist}(x,\partial\Omega)^2 \ \mbox{ uniformly on every compact subset of } \Omega,
\end{equation}
where the function $\Phi$ is defined by
\begin{equation}
\label{definition of Phi}
\Phi(s) = \int_1^s \frac {\phi^\prime(\xi)}{\xi} d\xi\quad\mbox{ for } s > 0.
\end{equation}
However, when $1 < p < 2$ or $0 < m <1$, the short time behavior of $u$ is not controlled by the distance function in such a way. To overcome this difficulty in the proofs of Theorems \ref{th:interaction curvatures pLaplace} and \ref{th:interaction curvatures mPorous},  we use the fact that the short time behavior of the solution $u$ is described by the boundary blow-up solutions given in \cite{Mjam1996, BMjam1992}. The results of the present paper in the case where $f\equiv g \equiv 1$ were announced in \cite{Soberwolfach2013}.

The present paper is organized as follows.  Section \ref{section1.5} is devoted to some preliminaries; the definitions of bounded solutions are mentioned,  the regularity results for the solutions  are quoted from the references,  and we refer to the references for the comparison principles.  Throughout the following four sections the comparison principles, which are mentioned in Section \ref{section1.5},  play a key role.  In Section \ref{section2},  it is shown that the short time behavior of the solutions is described by the boundary blow-up solutions given in \cite{Mjam1996, BMjam1992} in the case where  $\partial\Omega$ is bounded and of class $C^2$ and where both $f$ and $g$ are positive constants. 
In Section \ref{section3},  the problems are reduced to the case where $\partial\Omega$ is bounded and of class $C^2$ and where both $f$ and $g$ are positive constants.  Sections \ref{section4} and \ref{section4.5} are devoted to the construction of super- and subsolutions near the boundary $\partial\Omega$ for short times in the $p$-Laplace case and  in the porous medium type case, respectively.  
  In Section \ref{section5} we prove Theorems \ref{th:interaction curvatures pLaplace} and \ref{th:interaction curvatures mPorous}.


\setcounter{equation}{0}
\setcounter{theorem}{0}

\section{Prelimiaries: bounded solutions, regularity and comparison principles}
\label{section1.5}

Let us first consider the equation $\partial_t u= \mbox{\rm div}(|\nabla u|^{p-2}{\nabla u})$ with $1<p<2$.
By a bounded solution $u$ of problem {\rm (\ref{diffusion pLaplace})-(\ref{initial pLaplace})} we mean that $u \in C^{0}(\overline{\Omega} \times (0,\infty)) \cap L^{p}_{loc}\left(0,\infty; W^{1,p}_{loc}(\Omega)\right) \cap L^{\infty}(\Omega\times(0,\infty))$ satisfies \eqref{diffusion pLaplace} in the weak sense and $u(\cdot,t) \to 0$ in $L^{1}_{loc}(\Omega)$ as $t \to 0^{+}$, and by a bounded solution $u$ of problem {\rm (\ref{cauchy pLaplace})} we mean that $u \in C_{loc}(0,\infty; L^{2}_{loc}(\mathbb R^{N})) \cap L^{p}_{loc}\left(0,\infty; W^{1,p}_{loc}(\mathbb R^{N})\right) \cap L^{\infty}(\mathbb R^{N}\times(0,\infty))$ satisfies the differential equation in the weak sense and $u(\cdot,t) \to g(\cdot)\mathcal X_{\Omega^c}(\cdot)$ in $L^{1}_{loc}(\mathbb R^{N})$ as $t \to 0^{+}$. 

It is known that such bounded solutions $u$ together with $\nabla u$ are locally H\"older continuous, and both boundary and initial regularity of such solutions are known. See {\rm \cite{DiBspringer1993, DiBGV2012, Lnonlanal1993}}. Moreover, it is shown in {\rm \cite[Corollary 2.1, p. 2159]{BIVadvMath2010}} that such solutions are local strong ones, more precisely $\partial_{t} u \in L^{2}_{loc}$. 

The comparison principle for such strong solutions is obtained by Kurta {\rm \cite{KCRAcadSciParis1996, Kapplanal1999}}
for both the initial-boundary value problem and the Cauchy problem.  Furthermore,  note that one can easily prove Kurta's comparison principle  also for bounded weak solutions by taking his testing function modulo a Steklov time averaging process. 
See {\rm \cite{DiBspringer1993, DiBGV2012}} for the process, and 
see also {\rm \cite[Corollary 1.1, p. 189]{DiBGV2012}} for the comparison principle for weak solutions of the initial-boundary value problem over bounded domains.

Let us next consider the porous medium type equation $\partial_t u= \Delta u^{m}$ with $0<m<1$. 
By a bounded nonnegative solution $u$ of problem {\rm (\ref{diffusion mPorous})-(\ref{initial mPorous})} we mean that $u \in C^{0}(\overline{\Omega} \times (0,\infty)) \cap L^{\infty}(\Omega\times(0,\infty))$ is nonnegative and satisfies \eqref{diffusion mPorous} in the weak sense and $u(\cdot,t) \to 0$ in $L^{1}_{loc}(\Omega)$ as $t \to 0^{+}$, and by a bounded nonnegative solution $u$ of problem {\rm (\ref{cauchy mPorous})} we mean that $u \in C_{loc}(0,\infty; L^{2}_{loc}(\mathbb R^{N})) \cap L^{\infty}(\mathbb R^{N}\times(0,\infty))$ is nonnegative and satisfies the differential equation in the weak sense and $u(\cdot,t) \to g(\cdot)\mathcal X_{\Omega^c}(\cdot)$ in $L^{1}_{loc}(\mathbb R^{N})$ as $t \to 0^{+}$. 

It is known that such bounded solutions $u$  are locally H\"older continuous, and both boundary and initial regularity of such solutions are known. See {\rm \cite{DiBspringer1993, DiBGV2012}}. 

The comparison principle for such solutions of both the initial-boundary value problem and the Cauchy problem can be easily proved by modifying the proofs of {\rm \cite[Theorem A.1, pp. 253--257]{MSjde2012}} and {\rm \cite[Proposition A, pp. 1006--1008]{BKPnonlanal1985}}, with the aid of an idea of Dahlberg and Kenig {\rm \cite[Lemma 2.3, pp. 271--273]{DaK1986}} which circumvents the singularity coming from $u^{m}$ with $0 < m < 1$ at $u=0$.
See also {\rm \cite[Corollary 5.1, p. 201]{DiBGV2012}} for the comparison principle for weak solutions of the initial-boundary value problem over bounded domains.


\setcounter{equation}{0}
\setcounter{theorem}{0}

\section{Initial behavior and boundary blow-up solutions}
\label{section2}

Let $\Omega$ be a domain in $\mathbb R^N$ where $ \partial\Omega$  is bounded and of class $C^2$. Then it is known that there  exists a unique solution $v \in W^{1,p}_{loc}(\Omega)$ of
\begin{eqnarray}
&\mbox{\rm div}(|\nabla v|^{p-2}{\nabla v}) = \frac 1{2-p} v \ &\mbox{ and }\ v > 0\ \mbox{ in }\ \Omega,\label{pLaplace}\\
&v(x) \to \infty\ \ &\mbox{as  }\ x \to \partial\Omega,\label{boundary blow-up pLaplace}\\
&v(x) \to 0\ \ &\mbox{as }\ |x| \to \infty\  \mbox{ provided $\Omega$ is unbounded.}\label{at infinity pLaplace}
\end{eqnarray}
Here, $v$ belongs to $C^{1}(\Omega)$ and $\nabla v$ is locally H\"older continuous in $\Omega$, and moreover
\begin{equation}
\label{boundary behavior for pLaplace}
\frac {v(x)}{d(x)^{-\frac p{2-p}}} \to c(p)\ \mbox{ as }\ d(x) \to 0\ \mbox{ uniformly in } \Omega,
\end{equation}
where 
\begin{equation}
\label{distance and c(p)}
d(x) = \mbox{ dist}(x,\partial\Omega) \mbox{ for } x \in \Omega\ \mbox{ and }\ c(p) = \frac {2-p}{p}\left(\frac{2-p}{2p(p-1)}\right)^{-\frac1{2-p}}.
\end{equation}
The case where $\Omega$ is bounded was proved in \cite[Theorem 6.4 and Corollary 4.5, p. 245 and p. 231]{Mjam1996} and the case where $\Omega$ is unbounded, that is, $\Omega$ is an exterior domain, the existence of $v$ can be obtained with the aid of the argument in \cite[1.6, p. 12]{BMjam1992}, and the uniqueness also follows  by virtue of  \eqref{at infinity pLaplace}. 

Also, it is known by \cite[Theorem 2.7, pp. 18--19]{BMjam1992} that there  exists a unique solution $w \in C^2(\Omega)$ of
\begin{eqnarray}
&\Delta w^m = \frac 1{1-m} w \ &\mbox{ and }\ w > 0\ \mbox{ in }\ \Omega,\label{mPorous}\\
&w(x) \to \infty\ \ &\mbox{as  }\ x \to \partial\Omega,\label{boundary blow-up mPorous}\\
&w(x) \to 0\ \ &\mbox{as }\ |x| \to \infty\ \mbox{ provided $\Omega$ is unbounded.}\label{at infinity mPorous}
\end{eqnarray}
Note that in \cite{BMjam1992} the function $w(x)^m$ is dealt with instead of $w(x)$. Moreover, 
\begin{equation}
\label{boundary behavior for mPorous}
\frac {w(x)}{d(x)^{-\frac 2{1-m}}} \to c(m)\ \mbox{ as }\ d(x) \to 0\ \mbox{ uniformly in } \Omega,
\end{equation}
where 
\begin{equation}
\label{c(m)}
 c(m) = \left(\frac{2m(1+m)}{1-m}\right)^{\frac1{1-m}}.
\end{equation}
See \cite[Theorem 2.3, p. 17]{BMjam1992} or \cite[Corollary 4.5, p. 231]{Mjam1996} for \eqref{boundary behavior for mPorous}. 


\begin{proposition}
\label{prop:initial behavior for pLaplace}
Assume that $\partial\Omega$  is bounded and of class $C^2$. Let $u$ be the solution of 
either problem {\rm (\ref{diffusion pLaplace})-(\ref{initial pLaplace})} or problem {\rm (\ref{cauchy pLaplace})} where both $f$ and $g$ are positive constants.
Then 
\begin{equation}
\label{initial behavior for pLaplace}
t^{-\frac{1}{2-p}}u(x,t) \to v(x)\ \mbox{ as } t \to 0^{+}\ \mbox{ uniformly on compact sets in }\Omega,
\end{equation}
and moreover
\begin{equation}
\label{upper bound by separable solution for pLaplace}
u(x,t) \le t^{\frac{1}{2-p}} v(x)\ \mbox{ in } \Omega \times (0,\infty),
\end{equation}
where $v$ is the solution of problem {\rm (\ref{pLaplace})-(\ref{at infinity pLaplace})}.
\end{proposition}

\noindent
{\it Proof.}\ Define the function $V = V(x,t)$ for $(x,t) \in \Omega \times (0,\infty)$ by
$$
V(x,t) = t^{\frac{1}{2-p}} v(x).
$$
Then $V$ solves 
\begin{eqnarray}
&V > 0\ \mbox{ and }\ \partial_t V=\mbox{\rm div}(|\nabla V|^{p-2}{\nabla V})\ \ &\mbox{in }\ \Omega\times (0,\infty),\label{diffusion pLaplace for V}
\\
&V=\infty\ \ &\mbox{on }\ \partial\Omega\times (0,\infty),\label{dirichlet pLaplace infty}
\end{eqnarray}
Therefore it follows from the comparison principle that
$$
u \le V\ \mbox{ in } \Omega \times (0,\infty),
$$
which gives \eqref{upper bound by separable solution for pLaplace}. 

Since $\partial\Omega$ is bounded and of class $C^2$, there exists a number $\varepsilon_{0}>0$ such that,
 for each $\varepsilon \in (0, \varepsilon_{0}]$, the set $\Omega^{\varepsilon}$ defined by
 \begin{equation}
 \label{little expanded domain}
 \Omega^{\varepsilon} = \{ x \in \mathbb R^{N}\ :\ \mbox{ dist}(x, \overline{\Omega}) < \varepsilon \}
 \end{equation}
 is also a domain with bounded $C^{2}$ boundary $\partial\Omega^{\varepsilon}$.  To distinguish the notation $\Omega^{\varepsilon}$ from the complement $\Omega^c= \mathbb R^N\setminus\Omega$,
 hereafter we never use the letter `` $c$ '' for this definition \eqref{little expanded domain}. 
 For each $\varepsilon \in (0, \varepsilon_{0})$, consider the boundary blow-up solution $v_{\varepsilon} \in  C^{1}(\Omega^{\varepsilon})$ of
\begin{eqnarray}
&\mbox{\rm div}(|\nabla v_{\varepsilon}|^{p-2}{\nabla v_{\varepsilon}}) = \frac 1{2-p} v_{\varepsilon} \ &\mbox{ and }\ v_{\varepsilon} > 0\ \mbox{ in }\ \Omega^{\varepsilon},\label{pLaplace epsilon epsilon}\\
&v_{\varepsilon}(x) \to \infty\ \ &\mbox{as  }\ x \to \partial\Omega^{\varepsilon},\label{boundary blow-up pLaplace epsilon}\\
&v_{\varepsilon}(x) \to 0\ \ &\mbox{as }\ |x| \to \infty\  \mbox{ provided $\Omega^{\varepsilon}$ is unbounded.}\label{at infinity pLaplace epsilon}
\end{eqnarray}
In view of the argument in \cite[Proof of Theorem 4.4, pp. 239--240]{Mjam1996}, we observe that $\mbox{ dist}(x,\partial\Omega^{\varepsilon}) = \mbox{ dist}(x,\partial\Omega) + \varepsilon$ for $x \in \Omega$ and there exists $r > 0$ independent of $\varepsilon$ such that $\Omega^{\varepsilon}$ satisfies the uniform interior and exterior ball condition with radius $r$ for all $\varepsilon \in (0, \varepsilon_{0}]$, and we see that 
\begin{equation}
\label{uniform approximation by distance epsilon}
v_{\varepsilon}  \to v\ \mbox{ as } \varepsilon \to 0^{+}\ \mbox{ uniformly on compact sets in } \Omega.
\end{equation}
Define the function $V_{\varepsilon} = V_{\varepsilon}(x,t)$ for $(x,t) \in \Omega^{\varepsilon} \times (0,\infty)$ by
$$
V_{\varepsilon}(x,t) = t^{\frac{1}{2-p}} v_{\varepsilon}(x).
$$
Then, for each $\varepsilon \in (0, \varepsilon_{0}]$, $V_{\varepsilon}$ solves 
\begin{eqnarray}
&V_{\varepsilon} > 0\ \mbox{ and }\ \partial_t V_{\varepsilon}=\mbox{\rm div}(|\nabla V_{\varepsilon}|^{p-2}{\nabla V_{\varepsilon}})\ \ &\mbox{in }\ \Omega^{\varepsilon}\times (0,\infty),\label{diffusion pLaplace for V epsilon}
\\
&V_{\varepsilon}=\infty\ \ &\mbox{on }\ \partial\Omega^{\varepsilon}\times (0,\infty),\label{dirichlet pLaplace infty epsilon}
\\
& V_{\varepsilon} = 0 \ \ &\mbox{in }\ \overline{\Omega^{\varepsilon/2}}. \label{initial pLaplace infty epsilon}
\end{eqnarray}
Hence, for each $\varepsilon \in (0, \varepsilon_{0}]$, there exists $t_{\varepsilon} > 0$ such that
$$
V_{\varepsilon} \le u \ \left\{\begin{array}{rl}
&\mbox{on } \partial\Omega \times (0,t_{\varepsilon}]\ \mbox{ if $u$ solves problem {\rm (\ref{diffusion pLaplace})-(\ref{initial pLaplace})}},
\\
&\mbox{on } \partial\Omega^{\varepsilon/2} \times (0,t_{\varepsilon}]\ \mbox{ if $u$ solves problem {\rm (\ref{cauchy pLaplace})},}
\end{array}\right.
$$
since both $f$ and $g$ are positive constants and both $\partial\Omega$ and $\partial\Omega^{\varepsilon/2}$
are compact sets in $\Omega^{\varepsilon}$. Thus, we have from the comparison principle 
$$
V_{\varepsilon} \le u \ \mbox{ in } \Omega \times (0, t_{\varepsilon}],
$$
which together with \eqref{upper bound by separable solution for pLaplace} concludes that
$$
v_{\varepsilon}(x) \le t^{-\frac{1}{2-p}}u(x,t) \le v(x)\ \mbox{ for every } (x,t) \in \Omega \times (0, t_{\varepsilon}].
$$
Therefore \eqref{initial behavior for pLaplace} follows from \eqref{uniform approximation by distance epsilon}. \qed


\begin{proposition}
\label{prop:initial behavior for mPorous}
Assume that $\partial\Omega$  is bounded and of class $C^2$. Let $u$ be the solution of 
either problem {\rm (\ref{diffusion mPorous})-(\ref{initial mPorous})} or problem {\rm (\ref{cauchy mPorous})} where both $f$ and $g$ are positive constants.
Then 
\begin{equation}
\label{initial behavior for mPorous}
t^{-\frac{1}{1-m}}u(x,t) \to w(x)\ \mbox{ as } t \to 0^{+}\ \mbox{ uniformly on compact sets in }\Omega,
\end{equation}
and moreover
\begin{equation}
\label{upper bound by separable solution for mPorous}
u(x,t) \le t^{\frac{1}{1-m}} w(x)\ \mbox{ in } \Omega \times (0,\infty),
\end{equation}
where $w$ is the solution of problem {\rm (\ref{mPorous})-(\ref{at infinity mPorous})}.
\end{proposition}

\noindent
{\it Proof.}\ This follows from the same argument as in the proof of Proposition {\rm \ref{prop:initial behavior for pLaplace}}. \qed


\setcounter{equation}{0}
\setcounter{theorem}{0}

\section{Reduction to the case where $\partial\Omega$ is bounded and of class $C^2$ and where both $f$ and $g$ are positive constants}
\label{section3}

Let us first consider the solution $u$ of 
 problem {\rm (\ref{diffusion pLaplace})-(\ref{initial pLaplace})}.
Let  $\alpha >  \frac {(N+1)(2-p)}{2p},\ x_0 \in \Omega$, and  assume that $B_R(x_0)$ is contained in  $\Omega$ and such that $\overline{B_R(x_0)} \cap \partial\Omega = \{ y_0 \}$ for some $y_0 \in \partial\Omega$ and $\partial\Omega \cap B_{\delta}(y_0)$ is of class $C^{2}$ for some $\delta >0$. We find a bounded $C^{2}$ domain
$\Omega_*$ satisfying
$$
B_{R}(x_{0}) \subset \Omega_* \subset \Omega,\ \overline{B_R(x_0)} \cap \partial\Omega_* = \{ y_0 \},  \mbox{ and } B_{\frac \delta2}(y_{0}) \cap \partial\Omega \subset \partial\Omega_*\cap \partial\Omega \subset B_{\delta}(y_{0}) \cap \partial\Omega.
$$
Let $\hat{u}=\hat{u}(x,t)$ be  the bounded solution of 
  the initial-boundary value problem:
  \begin{eqnarray}
&\partial_t \hat{u}=\mbox{\rm div}(|\nabla \hat{u}|^{p-2}{\nabla \hat{u}})\ \ &\mbox{in }\ \Omega_*\times (0,\infty),\label{diffusion pLaplace on hatOmega}\\
&\hat{u}=\max f\ \ &\mbox{on }\ \partial\Omega_*\times (0,\infty),\label{dirichlet pLaplace on hatOmega}\\
&\hat{u}=0\ \ &\mbox{on }\ \Omega_*\times \{0\}.\label{initial pLaplace on hatOmega}
\end{eqnarray}
Then by the comparison principle we have
\begin{equation}
\label{upper bound on hatOmega}
u \le \hat{u}\ \mbox{ in } \Omega_*\times (0,\infty).
\end{equation}

Take a small $\varepsilon > 0$ arbitrarily. 
Choose a function $\hat{f_{\varepsilon}} \in C^{2}(\partial\Omega_*)$ satisfying
\begin{equation}
\label{a function hat{f}epsilon}
\hat{f_{\varepsilon}}(y_{0}) = f(y_{0}) + \frac \varepsilon2,\ \hat{f_{\varepsilon}}=\max f+ \frac \varepsilon2\ \mbox{ on } \Omega \cap \partial\Omega_*, \mbox{ and } \hat{f_{\varepsilon}} \ge f\ \mbox{ on } \partial\Omega \cap \partial\Omega_*.
\end{equation}
Let $\hat{v_{\varepsilon}} \in C^{1}(\overline{\Omega_*})$ solve
 \begin{eqnarray}
&0=\mbox{\rm div}(|\nabla \hat{v_{\varepsilon}}|^{p-2}{\nabla \hat{v_{\varepsilon}}})\ \ &\mbox{in }\ \Omega_*,\label{pLaplace on hatOmega}\\
&\hat{v_{\varepsilon}}=\hat{f_{\varepsilon}}\ \ &\mbox{on }\ \partial\Omega_*.\label{pLaplace Dirichlet on hatOmega}
\end{eqnarray}
Then by the comparison principle we have
\begin{equation}
\label{stationary upper bound on hatOmega}
u \le \hat{v_{\varepsilon}}\ \mbox{ in } \Omega_*\times (0,\infty).
\end{equation}
Moreover, we can find a small number $\delta_{\varepsilon}  \in (0,\delta/9)$ and two  $C^2$ domains $\Omega_{+,\varepsilon}$ and $\Omega_{-,\varepsilon}$ having bounded $C^{2}$ boundaries with the following properties:
both $\Omega_{+,\varepsilon}$ and $\mathbb R^N \setminus \overline{\Omega_{-,\varepsilon}}$ are bounded;  $\mathbb R^N \setminus \overline{\Omega_{-,\varepsilon}} \subset B_{3\delta_{\varepsilon}}(y_0);$ 
$B_R(x_0)\subset\Omega_{+,\varepsilon} \subset  \Omega_* \subset \Omega \subset \Omega_{-,\varepsilon};$
$\overline{B_R(x_0)} \cap \partial\Omega_{+,\varepsilon} = \overline{B_R(x_0)} \cap \partial\Omega_{-,\varepsilon}= \{ y_0 \};$ 
$\partial \Omega_{+,\varepsilon} \cap \partial\Omega_* \subset B_{2\delta_{\varepsilon}}(y_0) \cap \partial\Omega;$
\begin{equation*}
\label{two constructed domains}
B_{\delta_{\varepsilon}}(y_0) \cap \partial\Omega \subset \partial\Omega_{\pm,\varepsilon}\cap\partial\Omega
\subset B_{2\delta_{\varepsilon}}(y_0) \cap \partial\Omega \left(\subset \partial\Omega_*\cap \partial\Omega\right); \ 
 \end{equation*}
 \begin{equation}
 \label{sufficient estimates}
 f(y_0)-\varepsilon \le f \ \mbox{ on } \overline{B_{4\delta_{\varepsilon}}(y_0)}\cap \partial\Omega \ \mbox{ and }\ \hat{v_{\varepsilon}} \le f(y_0) + \varepsilon \mbox{ on } \overline{B_{4\delta_{\varepsilon}}(y_0)}\cap \overline{\Omega_*}.
 \end{equation}
 \par
Let $u_{\pm}^{\varepsilon}=u_{\pm}^{\varepsilon}(x,t)$ be the two bounded solutions of the initial-boundary value problems: 
\begin{eqnarray}
&\partial_t u_{\pm}^{\varepsilon}=\mbox{\rm div}(|\nabla u_{\pm}^{\varepsilon}|^{p-2}{\nabla u_{\pm}^{\varepsilon}})\ \ &\mbox{in }\ \Omega_{\pm,\varepsilon}\times (0,\infty),\label{diffusion pLaplace on Omega-pm}\\
&u_{\pm}^{\varepsilon}= f(y_{0}) \pm\varepsilon\ \ &\mbox{on }\ \partial{\Omega_{\pm,\varepsilon}}\times (0,\infty),\label{dirichlet pLaplace on Omega-pm}\\
&u_{\pm}^{\varepsilon}=0\ \ &\mbox{on }\ \Omega_{\pm,\varepsilon}\times \{0\}.\label{initial pLaplace on Omega-pm}
\end{eqnarray}
Here we obtain 
\begin{proposition}
\label{prop:estimates from above and below by pm-epsion}
Let  $u$ be the solution of 
 problem {\rm (\ref{diffusion pLaplace})-(\ref{initial pLaplace})}.
For every small $\varepsilon > 0$ there exists $\tau_{\varepsilon}>0$ satisfying
$$
u_{-}^{\varepsilon} \le u \le u_{+}^{\varepsilon} \ \mbox{ in } B_{R}(x_{0}) \times (0, \tau_{\varepsilon}],
$$
where $u_{\pm}^{\varepsilon}$ are the solutions of problems  {\rm (\ref{diffusion pLaplace on Omega-pm})-(\ref{initial pLaplace on Omega-pm})}.
\end{proposition}

\noindent
{\it Proof.}\ By combining \eqref{stationary upper bound on hatOmega} and the second inequality of \eqref{sufficient estimates} with \eqref{dirichlet pLaplace on Omega-pm}, we see that
\begin{equation}
\label{boundary comparison-1}
u \le u_{+}^{\varepsilon}\  \mbox{ on } \left(\partial \Omega_{+,\varepsilon} \cap \overline{B_{4\delta_{\varepsilon}}(y_0)}\right) \times (0, \infty).
\end{equation}
Since $\partial \Omega_{+,\varepsilon} \setminus B_{4\delta_{\varepsilon}}(y_0)$ is a compact set contained in $\Omega_*$,  by applying Proposition \ref{prop:initial behavior for pLaplace} to the bounded $C^2$ domain $\Omega_*$ and the solution $\hat{u}$ of problem {\rm (\ref{diffusion pLaplace on hatOmega})-(\ref{initial pLaplace on hatOmega})}, we have from the corresponding estimate \eqref{upper bound by separable solution for pLaplace} and \eqref{upper bound on hatOmega} that there exists $\tau_{1,\varepsilon} > 0$ satisfying
\begin{equation}
\label{estimate from above in short time1}
u \le u_{+}^{\varepsilon}\  \mbox{ on } \Bigl(\partial \Omega_{+,\varepsilon} \setminus B_{4\delta_{\varepsilon}}(y_0)\Bigr)\times (0, \tau_{1,\varepsilon}].
\end{equation}
Hence with the aid of \eqref{boundary comparison-1} and \eqref{estimate from above in short time1} we have from the comparison principle that 
\begin{equation}
\label{estimate from above in short time all}
u \le u_{+}^{\varepsilon} \ \mbox{ in }  \Omega_{+,\varepsilon} \times (0, \tau_{1,\varepsilon}].
\end{equation}

On the other hand, the first inequality of \eqref{sufficient estimates} gives
\begin{equation}
\label{boundary comparison-2}
u_{-}^{\varepsilon} \le u \ \mbox{ on } \left( \partial\Omega\cap\overline{B_{4\delta_{\varepsilon}}(y_0)}\right) \times (0,\infty).
\end{equation}
Since $\partial\Omega_{-,\varepsilon} \subset \overline{B_{3\delta_{\varepsilon}}(y_0)}$, by applying Proposition \ref{prop:initial behavior for pLaplace} to the domain $\Omega_{-,\varepsilon}$ with bounded $C^2$ boundary and the solution $u_{-}^{\varepsilon}$ of problem {\rm (\ref{diffusion pLaplace on Omega-pm})-(\ref{initial pLaplace on Omega-pm})}, we have from the corresponding estimate \eqref{upper bound by separable solution for pLaplace} and \eqref{f is bounded from above and below} that there exists $\tau_{2,\varepsilon} > 0$ satisfying
\begin{equation}
\label{estimate from below in short time2}
u_{-}^{\varepsilon} \le u\  \mbox{ on }  \left( \partial\Omega\setminus\overline{B_{4\delta_{\varepsilon}}(y_0)}\right) \times (0, \tau_{2,\varepsilon}].
\end{equation}
Therefore with the aid of \eqref{boundary comparison-2} and \eqref{estimate from below in short time2} we have from the comparison principle that 
\begin{equation}
\label{estimate from below in short time all}
u_{-}^{\varepsilon} \le u\  \mbox{ in }  \Omega \times (0, \tau_{2,\varepsilon}].
\end{equation}
In conclusion, \eqref{estimate from above in short time all} and \eqref{estimate from below in short time all} complete the proof if we set $\tau_\varepsilon = \min\left\{ \tau_{1,\varepsilon},\tau_{2,\varepsilon}\right\}$.
\qed

\vskip 5ex
Let us next consider the solution $u$ of  problem  {\rm (\ref{cauchy pLaplace})}.  Take a small $\varepsilon > 0$ arbitrarily.
Since $g(y_0) > 0$ and $g \in C^0(\mathbb R^N)$, there exists a small number $\delta_{\varepsilon}  \in (0,\delta/9)$ such that
\begin{equation}
\label{small nbd of yo}
g(y_0) - \frac 12\varepsilon \le g \le g(y_0) + \frac 12\varepsilon\ \mbox{ in } \overline{B_{4\delta_{\varepsilon}}(y_0)}.
\end{equation}
Moreover we  find a small number $\gamma_\varepsilon \in (0,\delta_\varepsilon)$ and two  $C^2$ domains $\Omega_{+,\varepsilon}$ and $\Omega_{-,\varepsilon}$ having bounded $C^{2}$ boundaries with the following properties:
both $\Omega_{+,\varepsilon}$ and $\mathbb R^N \setminus \overline{\Omega_{-,\varepsilon}}$ are bounded;  $\mathbb R^N \setminus \overline{\Omega_{-,\varepsilon}} \subset B_{3\delta_{\varepsilon}}(y_0);$ 
$B_R(x_0)\subset\Omega_{+,\varepsilon} \subset  \Omega \subset \Omega_{-,\varepsilon};$
$\overline{B_R(x_0)} \cap \partial\Omega_{+,\varepsilon} = \overline{B_R(x_0)} \cap \partial\Omega_{-,\varepsilon}= \{ y_0 \};$ 
$B_{\delta_{\varepsilon}}(y_0) \cap \partial\Omega \subset \partial\Omega_{\pm,\varepsilon}\cap\partial\Omega
\subset B_{2\delta_{\varepsilon}}(y_0) \cap \partial\Omega;$ 
\begin{equation}
\label{estimate from above can be by  local annulus}
 \overline{\left(\Omega_{+,\varepsilon}\right)^{\gamma_\varepsilon} } \cap \left(\mathbb R^N \setminus \Omega\right) \subset B_{4\delta_{\varepsilon}}(y_0),
 \end{equation}
where $\left(\Omega_{+,\varepsilon}\right)^{\gamma_\varepsilon}$ is the domain defined by \eqref{little expanded domain}, that is, 
$$
\left(\Omega_{+,\varepsilon}\right)^{\gamma_\varepsilon} =  \{ x \in \mathbb R^{N}\ :\ \mbox{ dist}(x, \overline{\Omega_{+,\varepsilon}}) < \gamma_\varepsilon \}.
$$
 \par
Let $u_{\pm}^{\varepsilon}=u_{\pm}^{\varepsilon}(x,t)$ be the two bounded solutions of the Cauchy problems \eqref{cauchy pLaplace} where the initial data $g\mathcal X_{\Omega^c}$ is 
replaced by $(g(y_0) \pm \varepsilon )\mathcal X_{\left(\Omega_{\pm,\varepsilon}\right)^c}$, respectively.
Hence we have
\begin{proposition}
\label{prop:estimates from above and below by pm-epsion cauchy}
Let  $u$ be the solution of 
 problem  {\rm (\ref{cauchy pLaplace})}.
For every small $\varepsilon > 0$ there exists $\tau_{\varepsilon}>0$ satisfying
$$
u_{-}^{\varepsilon} \le u \le u_{+}^{\varepsilon} \ \mbox{ in } B_{R}(x_{0}) \times (0, \tau_{\varepsilon}],
$$
where $u_{\pm}^{\varepsilon}$ are the solutions of problems \eqref{cauchy pLaplace} where the initial data $g\mathcal X_{\Omega^c}$ is 
replaced by $(g(y_0) \pm \varepsilon )\mathcal X_{\left(\Omega_{\pm,\varepsilon}\right)^c}$, respectively.
\end{proposition}

\noindent
{\it Proof.}\ In view of  \eqref{small nbd of yo} and the fact that $\mathbb R^N \setminus \overline{\Omega_{-,\varepsilon}} \subset B_{3\delta_{\varepsilon}}(y_0)$, we notice that
$$
(g(y_0) - \varepsilon )\mathcal X_{\left(\Omega_{-,\varepsilon}\right)^c} \le g\mathcal X_{\Omega^c}\ \mbox{ in } \mathbb R^N.
$$
Hence it follows from the comparison principle that
\begin{equation}
\label{cauchy from below by minus}
u_{-}^{\varepsilon} \le u\ \mbox{ in }\ \mathbb R^N \times (0,\infty).
\end{equation}

On the other hand, \eqref{small nbd of yo} and \eqref{estimate from above can be by  local annulus} yield that
$$
 g\mathcal X_{\Omega^c} \le g(y_0) + \frac 12\varepsilon < g(y_0) + \varepsilon = (g(y_0) + \varepsilon )\mathcal X_{\left(\Omega_{+,\varepsilon}\right)^c}\ \mbox{ in }
\left(\Omega_{+,\varepsilon}\right)^{\gamma_\varepsilon} \setminus \Omega_{+,\varepsilon}.
$$
Therefore by the initial behavior of the solutions there exists $\tau_\varepsilon > 0$ such that
$$
u \le u_{+}^{\varepsilon}\ \mbox{ on } \partial \left(\Omega_{+,\varepsilon}\right)^{\gamma_\varepsilon/2} \times (0, \tau_\varepsilon],
$$
which together with the comparison principle yields that
\begin{equation}
\label{cauchy from above by plus}
u \le u_{+}^{\varepsilon}\ \ \mbox{ in }\ \left(\Omega_{+,\varepsilon}\right)^{\gamma_\varepsilon/2} \times (0,\tau_\varepsilon].
\end{equation}
Thus, combining \eqref{cauchy from below by minus} with \eqref{cauchy from above by plus} completes the proof. \qed

\vskip 3ex 
Finally, Propositions \ref{prop:estimates from above and below by pm-epsion} and \ref{prop:estimates from above and below by pm-epsion cauchy} yield
$$
\!\!\!\int\limits_{B_R(x_0)}\! \left(u_{-}^{\varepsilon}(x,t)\right)^{\alpha}\ dx \le\!\!\!\int\limits_{B_R(x_0)}\! \left(u(x,t)\right)^{\alpha}\ dx \le \!\!\!\int\limits_{B_R(x_0)}\! \left(u_{+}^{\varepsilon}(x,t)\right)^{\alpha}\ dx\ \mbox{ for every } t \in (0, \tau_{\varepsilon}].
$$
These two inequalities show that the proofs of Theorem \ref{th:interaction curvatures pLaplace} for  the equation $\partial_t u= \mbox{\rm div}(|\nabla u|^{p-2}{\nabla u})$ are reduced to the case where $\partial\Omega$ is bounded and of class $C^2$ and where $f$ and $g$  are  positive constants, since we later know that the positive constants $c$ in formula \eqref{asymptotics and curvatures pLaplace} are continuous with respect to positive constants $f$ and $g$, respectively.
Also, the proofs for the equation $\partial_t u= \Delta u^{m}$ follow from the same arguments as in those for the equation $\partial_t u= \mbox{\rm div}(|\nabla u|^{p-2}{\nabla u})$.


\setcounter{equation}{0}
\setcounter{theorem}{0}

\section{Super- and  subsolutions near the boundary for short times: the $p$-Laplace case}
\label{section4}

By virtue of section \ref{section3}, we can assume that $\partial\Omega$ is bounded and of class $C^2$ and  $f \equiv g \equiv \beta$ for some positive constant $\beta > 0$.

Let us first consider the solution $u$ of 
 problem {\rm (\ref{diffusion pLaplace})-(\ref{initial pLaplace})}.  Namely, we consider the bounded solution $u=u(x,t)$  of 
  the  initial-boundary value problem:
  \begin{eqnarray*}
&\partial_t u=\mbox{\rm div}(|\nabla u|^{p-2}{\nabla u})\ \ &\mbox{in }\ \Omega\times (0,\infty),\label{diffusion pLaplace beta}\\
&u=\beta\ \ &\mbox{on }\ \partial\Omega\times (0,\infty),\label{dirichlet pLaplace beta}\\
&u=0\ \ &\mbox{on }\ \Omega\times \{0\}.\label{initial pLaplace beta}
\end{eqnarray*}
For $\xi \ge 0$, define $\varphi = \varphi(\xi)$ by
\begin{equation}
\label{1-dimensional similarity solution for pLaplace ib-problem}
\varphi(\xi) = \beta - \left(\frac {2-p}{2p(p-1)}\right)^{-\frac 1{2-p}} \int_0^\xi (\eta^2 + \lambda)^{-\frac 1{2-p}} d\eta,
\end{equation}
where $\lambda > 0$ is determined uniquely by the equation $\varphi(\infty) = 0$. Then $\varphi = \varphi(\xi)$ satisfies
\begin{eqnarray}
&&(p-1) |\varphi^\prime|^{p-2}\varphi^{\prime\prime} + \frac 1p\varphi^\prime\xi = 0\ \mbox{ for } \xi > 0, \label{2nd order ODE for pLaplace ib-problem}
\\
&& \varphi(0) = \beta,\ \varphi^\prime < 0\ \mbox{ in } [0,\infty), \mbox{ and } \varphi(\infty) = 0.\label{additional for pLaplace ib-problem}
\end{eqnarray}
l'Hospital's rule gives
\begin{equation}
\label{limits corresponds to boundary blow up sol}
\lim_{\xi \to \infty} \frac {\varphi(\xi)}{\xi^{-\frac p{2-p}}} = c(p),
\end{equation}
where $c(p)$ is the constant given by \eqref{distance and c(p)}. Note that, if we set $h(s,t) = \varphi(t^{- 1/p}s) $ for $s \ge 0$ and $t > 0$, then $h$ satisfies the one-dimensional problem:
$$
\partial_t h = \partial_s\left( |\partial_s h|^{p-2}\partial_s h\right) \ \mbox{ in } (0,\infty)^2,\ h = \beta\ \mbox{ on } \{0\} \times (0,\infty), \ \mbox{ and } h = 0\ \mbox{ on }(0,\infty)\times\{0\}.
$$
For small $\varepsilon >0$, define $\varphi_\pm = \varphi_\pm(\xi)\ (\xi > 0)$ by
\begin{equation}
\label{1-dimensional similarity solution for pLaplace ib-problem pm}
\varphi_\pm(\xi) = \beta \pm\varepsilon - \left(\frac {2-p}{2p(p-1)}\right)^{-\frac 1{2-p}} \int_0^\xi \left(\eta^2 \mp 2p\varepsilon\int_0^\eta\sqrt{1+s^2} ds+ \lambda_\pm\right)^{-\frac 1{2-p}} d\eta,
\end{equation}
where each $\lambda_\pm > 0$ is determined uniquely by the equation $\varphi_\pm(\infty) = 0.$ Notice that  
\begin{eqnarray}
&& \varphi_{\pm}\to \varphi\ \mbox{ as } \varepsilon \to 0^{+} \mbox{ uniformly on } [0,\infty),\label{approximation for varphi}
\\
&&\lambda_{\pm} \to \lambda\ \mbox{ as } \varepsilon \to 0^{+},\label{approximation for lambda} 
\end{eqnarray}
where $\lambda$ is given in \eqref{1-dimensional similarity solution for pLaplace ib-problem}.
Then $\varphi_\pm= \varphi_\pm(\xi)$ satisfies
\begin{eqnarray}
&&(p-1) |\varphi_\pm^\prime|^{p-2}\varphi_\pm^{\prime\prime} + \frac 1p\varphi_\pm^\prime\left[\xi \mp p\varepsilon\sqrt{1+\xi^2}\right]= 0\ \mbox{ for } \xi > 0, \label{2nd order ODE for pLaplace ib-problem pm}
\\
&& \varphi_\pm(0) = \beta\pm\varepsilon,\ \varphi_\pm^\prime < 0\ \mbox{ in } [0,\infty), \mbox{ and } \varphi_\pm(\infty) = 0.\label{additional for pLaplace ib-problem pm}
\end{eqnarray}
l'Hospital's rule gives
\begin{equation}
\label{limits corresponds to to boundary blow up sol pm}
\lim_{\xi \to \infty} \frac {\varphi_\pm(\xi)}{\xi^{-\frac p{2-p}}} = c(p)(1\mp p\varepsilon)^{-\frac1{2-p}}.
\end{equation}
Since $\partial\Omega$ is bounded and of class $C^2$, there exists $\rho_0 > 0$  such that the distance function $d = d(x)$ of $x \in \overline{\Omega}$ to the boundary $\partial\Omega$ is $C^2$-smooth on $\overline{\Omega_{\rho_0}}$, 
where  $\Omega_{\rho_0}$ is defined by  \eqref{neighborhood of boundary} with $\rho = \rho_0$.

By setting
\begin{equation}
\label{super&subsolutions for pLaplace ib-problem}
w_\pm(x,t) =\varphi_\pm(t^{-1/p}d(x))\ \mbox{ for } (x,t) \in \Omega \times (0,\infty),
\end{equation}
we obtain
\begin{proposition}
\label{prop:estimates from above and below nbd of boundary by pm}
Let  $u$ be the solution of 
 problem {\rm (\ref{diffusion pLaplace})-(\ref{initial pLaplace})} where $\partial\Omega$ is bounded and of class $C^2$ and  $f \equiv \beta$ for some positive constant $\beta > 0$. 
 For every small $\varepsilon > 0$ there exist $\rho_\varepsilon \in (0, \rho_0)$ and $\tau_{\varepsilon}>0$ satisfying
\begin{equation}
\label{important estimates for the solution near the boundary for short time}
w_{-} \le u \le w_{+} \ \mbox{ in } \Omega_{\rho_\varepsilon} \times (0, \tau_{\varepsilon}],
\end{equation}
where $w_{\pm}$ are given by \eqref{super&subsolutions for pLaplace ib-problem} and $\Omega_{\rho_\varepsilon}$ is defined by  \eqref{neighborhood of boundary} with $\rho = \rho_\varepsilon$.
\end{proposition}

\noindent
{\it Proof.}\ Take a small $\varepsilon > 0$. For $x \in \Omega_{\rho_0}$ and $t >0$, a straightforward computation gives
$$
\partial_t w_\pm - \mbox{\rm div}(|\nabla w_\pm|^{p-2}{\nabla w_\pm}) = - t^{-1}\varphi_\pm^\prime\left[ \pm\varepsilon\sqrt{1+ \xi^2} + t^{1/p} |\varphi_\pm^\prime |^{p-2}\Delta d\right],
$$
where $\xi = t^{-1/p} d(x)$ and 
$$
|\varphi_\pm^\prime|^{p-2} = (-\varphi_\pm^\prime)^{p-2} = \left(\frac{2-p}{2p(p-1)}\right) \left[\xi^2 \mp 2p\varepsilon\int_0^\xi\sqrt{1+s^2}ds + \lambda_\pm\right].
$$
Therefore, by using \eqref{approximation for lambda} and observing that
$$
t^{1/p}\xi^{2} \le |\xi|d(x)\ \mbox{ and }\ t^{1/p}\left|\int_0^\xi\sqrt{1+s^2}ds\right| \le t^{1/p}(|\xi| + \xi^{2}),
$$
we notice that there exist $\rho_{1,\varepsilon} \in (0, \rho_0)$ and $\tau_{1,\varepsilon}>0$ satisfying
\begin{equation}
\label{differential inequalities for pLaplace initial-boundary}
(\pm1)\left(\partial_t w_\pm - \mbox{\rm div}(|\nabla w_\pm|^{p-2}{\nabla w_\pm})\right) > 0
 \ \mbox{ in } \Omega_{\rho_{1,\varepsilon}} \times (0, \tau_{1,\varepsilon}],
\end{equation}
where $w_{\pm}$ are given by \eqref{super&subsolutions for pLaplace ib-problem} and $\Omega_{\rho_{1,\varepsilon}}$ is defined by  \eqref{neighborhood of boundary} with $\rho = \rho_{1,\varepsilon}$.

By \eqref{boundary behavior for pLaplace}, there exists $\rho_{\varepsilon} \in (0, \rho_{1,\varepsilon})$ satisfying
\begin{equation*}
c(p)\left(1+\frac {p\varepsilon}4\right)^{-\frac 1{2-p}} d(x)^{-\frac p{2-p}} \le v(x) \le 
c(p)\left(1-\frac {p\varepsilon}4\right)^{-\frac 1{2-p}} d(x)^{-\frac p{2-p}}\ \mbox{ for } x \in \Omega_{\rho_{\varepsilon}}.
\end{equation*}
Hence by \eqref{initial behavior for pLaplace} of Proposition \ref{prop:initial behavior for pLaplace}
there exists $\tau_{2,\varepsilon} \in (0, \tau_{1,\varepsilon}]$ such that for $(x,t) \in \Gamma_{\rho_{\varepsilon}}\times (0, \tau_{2,\varepsilon}]$
\begin{equation}
\label{1st inequality for initial-boundary pLaplcae}
c(p)\left(1+\frac {p\varepsilon}2\right)^{-\frac 1{2-p}} (\rho_{\varepsilon})^{-\frac p{2-p}} \le t^{-\frac 1{2-p}}u(x,t) \le 
c(p)\left(1-\frac {p\varepsilon}2\right)^{-\frac 1{2-p}} (\rho_{\varepsilon})^{-\frac p{2-p}},
\end{equation}
where $\Gamma_{\rho_{\varepsilon}}$ is defined by  \eqref{level surfaces of distance functions} with $\rho = \rho_{\varepsilon}$.

Moreover, by \eqref{limits corresponds to to boundary blow up sol pm}, there exists $\tau_{\varepsilon} \in (0, \tau_{2,\varepsilon}]$ such that for $(x,t) \in \Gamma_{\rho_{\varepsilon}}\times (0, \tau_{\varepsilon}]$
\begin{eqnarray*}
&& t^{-\frac 1{2-p}}(\rho_{\varepsilon})^{\frac p{2-p}} w_{+}(x,t) \ge c(p)\left(1-\frac {p\varepsilon}2\right)^{-\frac 1{2-p}},
\\
&&t^{-\frac 1{2-p}}(\rho_{\varepsilon})^{\frac p{2-p}} w_{-}(x,t) \le c(p)\left(1+\frac {p\varepsilon}2\right)^{-\frac 1{2-p}}.
\end{eqnarray*}
Thus combining these inequalities with \eqref{1st inequality for initial-boundary pLaplcae} yields that
\begin{equation}
\label{inequalities on Gamma_{epsilon}}
w_{-}\le u \le w_{+}\ \mbox{ on } \Gamma_{\rho_{\varepsilon}} \times (0,\tau_{\varepsilon}].
\end{equation}

Observe that
\begin{eqnarray}
&w_{-}=\beta-\varepsilon < \beta=u < \beta+\varepsilon = w_{+}\ &\mbox{ on }\partial\Omega \times (0,\tau_{\varepsilon}],\label{comparison on the boundary original}
\\
& w_{-}=u=w_{+}=0\ &\mbox{ on } \Omega_{\rho_{\varepsilon}} \times \{0\}.\label{initial condition comparison}
\end{eqnarray}
Therefore, by combining these with \eqref{inequalities on Gamma_{epsilon}} and \eqref{differential inequalities for pLaplace initial-boundary}, we get the conclusion \eqref{important estimates for the solution near the boundary for short time} from the comparison principle.
\qed

\vskip 4ex
Let us next consider the solution $u$ of 
 problem {\rm (\ref{cauchy pLaplace})}.   Namely, we consider the bounded solution $u=u(x,t)$  of 
  the Cauchy problem:
  \begin{equation*}
\label{cauchy pLaplace beta}
\partial_t u=\mbox{\rm div}(|\nabla u|^{p-2}{\nabla u})\ \mbox{ in }\ \mathbb R^N \times (0, \infty)\quad\mbox{ and }\ u = \beta\mathcal X_{\Omega^c}\ \mbox{ on }\ \mathbb R^N \times \{0\},
\end{equation*}
where  $\mathcal X_{\Omega^c}$ is the characteristic function of the set $\Omega^c = \mathbb R^N \setminus \Omega$.
For $\xi \in \mathbb R$, define $\psi = \psi(\xi)$ by
\begin{equation}
\label{1-dimensional similarity solution for pLaplace c-problem}
\psi(\xi) = \beta - \left(\frac {2-p}{2p(p-1)}\right)^{-\frac 1{2-p}} \int_{-\infty}^\xi (\eta^2 + \lambda)^{-\frac 1{2-p}} d\eta,
\end{equation}
where $\lambda > 0$ is determined uniquely by the equation $\psi(\infty) = 0$. Then $\psi = \psi(\xi)$ satisfies
\begin{eqnarray}
&&(p-1) |\psi^\prime|^{p-2}\psi^{\prime\prime} + \frac 1p\psi^\prime\xi = 0\ \mbox{ for } \xi \in \mathbb R, \label{2nd order ODE for pLaplace c-problem}
\\
&& \psi(-\infty) = \beta,\ \psi^\prime < 0\ \mbox{ in } \mathbb R, \mbox{ and } \psi(\infty) = 0.\label{additional for pLaplace c-problem}
\end{eqnarray}
l'Hospital's rule gives
\begin{equation}
\label{limits corresponds to boundary blow up sol for cauchy}
\lim_{\xi \to \infty} \frac {\psi(\xi)}{\xi^{-\frac p{2-p}}} = c(p),
\end{equation}
where $c(p)$ is the constant given by \eqref{distance and c(p)}. Note that, if we set $h(s,t) = \psi(t^{- 1/p}s) $ for $s \in \mathbb R$ and $t > 0$, then $h$ satisfies the one-dimensional problem:
$$
\partial_t h = \partial_s\left( |\partial_s h|^{p-2}\partial_s h\right) \ \mbox{ in } \mathbb R \times(0,\infty) \mbox{ and } h = \beta\mathcal X_{(-\infty,0]}\ \mbox{ on } \mathbb R\times\{0\}.
$$
For small $\varepsilon >0$, define $\psi_\pm = \psi_\pm(\xi)\ (\xi\in\mathbb R)$ by
\begin{equation}
\label{1-dimensional similarity solution for pLaplace c-problem pm}
\psi_\pm(\xi) = \beta \pm\varepsilon - \left(\frac {2-p}{2p(p-1)}\right)^{-\frac 1{2-p}} \int_{-\infty}^\xi \left(\eta^2 \mp 2p\varepsilon\int_0^\eta\sqrt{1+s^2} ds+ \lambda_\pm\right)^{-\frac 1{2-p}} d\eta,
\end{equation}
where each $\lambda_\pm > 0$ is determined uniquely by the equation $\psi_\pm(\infty) = 0.$ Notice that  
\begin{eqnarray}
&&\psi_{\pm} \to \psi\ \mbox{ as } \varepsilon \to 0^{+}\ \mbox{ uniformly on }\mathbb R,\label{approximation for psi}
\\
&&\lambda_{\pm} \to \lambda\ \mbox{ as } \varepsilon \to 0^{+}, \label{approximation for lambda for cauchy}
\end{eqnarray}
where $\lambda$ is given in \eqref{1-dimensional similarity solution for pLaplace c-problem}.
Then $\psi_\pm= \psi_\pm(\xi)$ satisfies
\begin{eqnarray}
&&(p-1) |\psi_\pm^\prime|^{p-2}\psi_\pm^{\prime\prime} + \frac 1p\psi_\pm^\prime\left[\xi \mp p\varepsilon\sqrt{1+\xi^2}\right]= 0\ \mbox{ for } \xi \in \mathbb R, \label{2nd order ODE for pLaplace c-problem pm}
\\
&& \psi_\pm(-\infty) = \beta\pm\varepsilon,\ \psi_\pm^\prime < 0\ \mbox{ in } \mathbb R, \mbox{ and } \psi_\pm(\infty) = 0.\label{additional for pLaplace c-problem pm}
\end{eqnarray}
l'Hospital's rule gives
\begin{equation}
\label{limits corresponds to to boundary blow up sol pm for cauchy}
\lim_{\xi \to \infty} \frac {\psi_\pm(\xi)}{\xi^{-\frac p{2-p}}} = c(p)(1\mp p\varepsilon)^{-\frac1{2-p}}.
\end{equation}
As in \cite{MSjde2012}, let us introduce the signed distance function $d^{*} = d^{*}(x)$ of $x \in \mathbb R^{N}$
to the boundary $\partial\Omega$ defined by
\begin{equation*}
\label{signed distance}
d^*(x) = \left\{\begin{array}{rll}
 \mbox{ dist}(x,\partial\Omega)\ &\mbox{ if }\ x \in \Omega,
\\
-\mbox{ dist}(x,\partial\Omega)\ &\mbox{ if  }\ x \not\in \Omega.
\end{array}\right.
\end{equation*}
For every $\rho > 0$, let $\mathcal N_{\rho}$ be a compact neighborhood of $\partial\Omega$ in $\mathbb R^{N}$ defined by
\begin{equation}
\label{neighborhood of boundary from both sides}
\mathcal N_{\rho} = \{ x \in \mathbb R^N : -\rho \le d^*(x) \le \rho \}.
\end{equation}
If $\partial\Omega$  is bounded and of class $C^2$,  there exists a number $\rho_0 > 0$ such that $d^*(x)$ is $C^2$-smooth on $\mathcal N_{\rho_{0}}$. For simplicity we have used the same letter $\rho_0 >0$ as in the previous case for problem {\rm (\ref{diffusion pLaplace})-(\ref{initial pLaplace})}.

By setting
\begin{equation}
\label{super&subsolutions for pLaplace c-problem}
w_\pm(x,t) =\psi_\pm(t^{-1/p}d^{*}(x))\ \mbox{ for } (x,t) \in \mathbb R^{N} \times (0,\infty),
\end{equation}
we obtain
\begin{proposition}
\label{prop:estimates from above and below nbd of boundary by pm for cauchy}
Let  $u$ be the solution of 
 problem {\rm (\ref{cauchy pLaplace})} where $\partial\Omega$ is bounded and of class $C^2$ and  $g \equiv \beta$ for some positive constant $\beta > 0$. 
 For every small $\varepsilon > 0$ there exist $\rho_\varepsilon \in (0, \rho_0)$ and $\tau_{\varepsilon}>0$ satisfying
\begin{equation}
\label{important estimates for the solution near the boundary for short time for cauchy}
w_{-} \le u \le w_{+} \ \mbox{ in } \mathcal N_{\rho_{\varepsilon}} \times (0, \tau_{\varepsilon}],
\end{equation}
where $w_{\pm}$ are given by \eqref{super&subsolutions for pLaplace c-problem} and $\mathcal N_{\rho_\varepsilon}$ is defined by  \eqref{neighborhood of boundary from both sides} with $\rho = \rho_\varepsilon$.
\end{proposition}

\noindent
{\it Proof.}\ The proof is similar to that of Proposition \ref{prop:estimates from above and below nbd of boundary by pm}. The ingredients \eqref{approximation for lambda}, \eqref{limits corresponds to to boundary blow up sol pm}, and \eqref{comparison on the boundary original}  are replaced by \eqref{approximation for lambda for cauchy}, \eqref{limits corresponds to to boundary blow up sol pm for cauchy}, and the corresponding inequalities on $\{ x \in \mathbb R^{N} : d^{*}(x) = -\rho_{\varepsilon}\} \times (0,\tau_{\varepsilon}]$,
respectively.
\qed


\setcounter{equation}{0}
\setcounter{theorem}{0}

\section{Super- and  subsolutions near the boundary for short times: the porous medium type case}
\label{section4.5}

By virtue of section \ref{section3}, we can assume that $\partial\Omega$ is bounded and of class $C^2$ and  $f \equiv g \equiv \beta$ for some positive constant $\beta > 0$.

Concerning $\partial_t u= \Delta u^{m}$ with $0<m<1$, the same constructions of super- and  subsolutions as in \cite{MSjde2012} work. Let $u=u(x,t)$ be the bounded solution of 
 problem {\rm \eqref{diffusion mPorous}-\eqref{initial mPorous}} where $f \equiv \beta$.
  Namely, we consider the bounded solution $u=u(x,t)$  of 
  the  initial-boundary value problem:
  \begin{eqnarray*}
&\partial_t u= \Delta u^{m}\ \ &\mbox{in }\ \Omega\times (0,\infty),\label{diffusion mPorous beta}\\
&u=\beta\ \ &\mbox{on }\ \partial\Omega\times (0,\infty),\label{dirichlet mPorous beta}\\
&u=0\ \ &\mbox{on }\ \Omega\times \{0\}.\label{initial mPorous beta}
\end{eqnarray*}
Let us set $\phi(s) = s^{m}$ for $s \ge 0$. We use a result from Atkinson and Peletier \cite{AParma1974}:  
for every $\gamma > 0$, there exists a unique $C^2$ solution $f_\gamma = f_\gamma(\xi)$ of the problem:
\begin{eqnarray}
&& \left( \phi^\prime(f_\gamma) f_\gamma^\prime\right)^\prime + \frac 12 \xi  f_\gamma^\prime = 0\ \mbox{ in }\ [0,\infty),
\label{ode selfsimilar}
\\
&& f_\gamma(0) = \gamma,\quad f_\gamma(\infty) = 0,\label{boundary conditions}
\\
&& f_\gamma^\prime < 0\ \mbox{ in }\ [0,\infty).
\label{monotonicity}
\end{eqnarray}
Moreover, \cite[Theorem 5 and its example 3, p. 388 and p. 390]{AParma1974} gives
\begin{equation}
\label{limits corresponds to boundary blow up sol for initial-boundary}
\lim_{\xi\to\infty}\frac {f_{\gamma}(\xi)}{\xi^{-\frac2{1-m}}} = c(m),
\end{equation}
where $c(m)$ is the constant given by \eqref{c(m)}. This behavior comes from the structure of the equation $\partial_t u= \Delta u^{m}$ with $0<m<1$, and it is different from that of the equation of the form
$\partial_{t}u=\Delta \phi(u)$ with $\delta_{1}\le \phi^{\prime}(s)\le \delta_{2}\ (s\in\mathbb R)$ for two positive constants $\delta_{1}, \delta_{2}$, which is treated in \cite[(3.15), p. 243]{MSjde2012}.
Note that, if we put $h(s,t) = f_\gamma\left(t^{-1/2}s\right)$ for $s \ge0$ and $t > 0$, then $h$ satisfies the one-dimensional problem:
$$
\partial_t h = \partial_s^2 \phi(h)\ \mbox{ in } (0,\infty)^2,\ h = \gamma\ \mbox{ on }\ \{0\} \times (0,\infty), \mbox{ and } h = 0\ \mbox{ on } (0,\infty)\times \{0\}.
$$
Let $0< \varepsilon < \frac 14$. Then, as in \cite[Proof of Lemma 3.1, pp. 242--244]{MSjde2012}, by continuity 
we can find a sufficiently small $0 < \eta_\varepsilon << \varepsilon$ and two $C^2$ functions $f_{\pm} = f_{\pm}(\xi)$ for $\xi \ge 0$ satisfying:
\begin{eqnarray*}
&& f_\pm(\xi) = f_{\beta\pm \varepsilon}\left(\sqrt{1\mp 2\eta_\varepsilon}\ \xi\right)\ \mbox{ if } \xi \ge \eta_\varepsilon; \label{scaling of f}
\\
&&f_\pm^\prime < 0\ \mbox{ in } [0,\infty);
\label{monotonicity both}
\\
&& f_- < f_\beta < f_+\ \mbox{ in } [0,\infty);
\label{from above and below}
\\
&& \left(\phi^\prime(f_\pm)f_\pm^\prime\right)^\prime + \frac 12\xi f_\pm^\prime = h_\pm(\xi) f_\pm^\prime\ \mbox{ in } [0,\infty),\label{ modified ode}
\end{eqnarray*}
where $h_\pm = h_\pm(\xi)$ are defined by
\begin{equation}
\label{definition of h_pm}
h_\pm(\xi) = \left\{\begin{array}{rl}   \pm \eta_\varepsilon\xi\ &\mbox{ if }\ \xi \ge \eta_\varepsilon,\\
  \pm\eta_\varepsilon^2\ &\mbox{ if }\ \xi \le \eta_\varepsilon.
\end{array}\right.
\end{equation}
(Here, in order to use the functions $h_\pm$ also for problem {\rm (\ref{cauchy mPorous})} later, we defined $h_\pm(\xi)$ for all $\xi \in \mathbb R$.)
The above construction of $f_\pm$ directly implies that
\begin{equation}
\label{approximation for f_{beta}}
f_{\pm} \to f_{\beta}\ \mbox{ as } \varepsilon \to 0^{+}\ \mbox{ uniformly on } [0,\infty).
\end{equation}
Moreover, by \eqref{limits corresponds to boundary blow up sol for initial-boundary} we have
\begin{equation}
\label{limits corresponds to boundary blow up sol pm for initial-boundary}
\lim_{\xi\to\infty} \frac {f_{\pm}(\xi)}{\xi^{-\frac 2{1-m}}} = c(m)(1\mp2\eta_{\varepsilon})^{-\frac1{1-m}}.
\end{equation}

By setting
\begin{equation}
\label{super&subsolutions for mPorous ib-problem}
w_\pm(x,t) =f_\pm(t^{-1/2}d(x))\ \mbox{ for } (x,t) \in \Omega \times (0,\infty),
\end{equation}
we obtain
\begin{proposition}
\label{prop:estimates from above and below nbd of boundary by pm for mPorous}
Let  $u$ be the solution of 
 problem {\rm (\ref{diffusion mPorous})-(\ref{initial mPorous})} where $\partial\Omega$ is bounded and of class $C^2$ and  $f \equiv \beta$ for some positive constant $\beta > 0$. 
 For every small $\varepsilon > 0$ there exist $\rho_\varepsilon \in (0, \rho_0)$ and $\tau_{\varepsilon}>0$ satisfying
\begin{equation}
\label{important estimates for the solution near the boundary for short time for mPorous ib}
w_{-} \le u \le w_{+} \ \mbox{ in } \Omega_{\rho_\varepsilon} \times (0, \tau_{\varepsilon}],
\end{equation}
where $w_{\pm}$ are given by \eqref{super&subsolutions for mPorous ib-problem} and $\Omega_{\rho_\varepsilon}$ is defined by  \eqref{neighborhood of boundary} with $\rho = \rho_\varepsilon$.
\end{proposition}

\noindent
{\it Proof.}\ Take a small $\varepsilon > 0$. For $x \in \Omega_{\rho_0}$ and $t >0$, a straightforward computation gives
$$
\partial_t w_\pm - \Delta (w_\pm)^{m} = - t^{-1}f_\pm^\prime\left[ h_{\pm}(\xi) + t^{1/2}m(f_{\pm})^{-(1-m)}\Delta d\right],
$$
where $\xi = t^{-1/2} d(x)$. In view of \eqref{limits corresponds to boundary blow up sol pm for initial-boundary}, we observe that there exists a constant $C_{\varepsilon} > 0$ satisfying
$$
t^{1/2}m(f_{\pm})^{-(1-m)} \le \left\{\begin{array}{rl}
& t^{1/2}C_{\varepsilon} \xi^{2} = C_{\varepsilon} \xi d(x) \quad \mbox{ if } \xi \ge \eta_{\varepsilon},
\\
& t^{1/2}m(f_{\pm}(\sqrt{1\mp2\eta_{\varepsilon}}\eta_{\varepsilon}))^{-(1-m)} \le t^{1/2}C_{\varepsilon} {\eta_{\varepsilon}}^{2} \quad \mbox{ if } \xi \le \eta_{\varepsilon}.
\end{array}\right.
$$
Therefore, with the aid of the definition \eqref{definition of h_pm} of $h_{\pm}(\xi)$, we notice that there exist $\rho_{1,\varepsilon} \in (0, \rho_0)$ and $\tau_{1,\varepsilon}>0$ satisfying
\begin{equation}
\label{differential inequalities for mPorous initial-boundary}
(\pm1)\left(\partial_t w_\pm - \Delta (w_\pm)^{m}\right) > 0
 \ \mbox{ in } \Omega_{\rho_{1,\varepsilon}} \times (0, \tau_{1,\varepsilon}],
\end{equation}
where $w_{\pm}$ are given by \eqref{super&subsolutions for mPorous ib-problem} and $\Omega_{\rho_{1,\varepsilon}}$ is defined by  \eqref{neighborhood of boundary} with $\rho = \rho_{1,\varepsilon}$.

By \eqref{boundary behavior for mPorous}, there exists $\rho_{\varepsilon} \in (0, \rho_{1,\varepsilon})$ satisfying
\begin{equation*}
c(m)\left(1+\frac {\eta_\varepsilon}2\right)^{-\frac 1{1-m}} d(x)^{-\frac 2{1-m}} \le w(x) \le 
c(m)\left(1-\frac {\eta_\varepsilon}2\right)^{-\frac 1{1-m}} d(x)^{-\frac 2{1-m}}\ \mbox{ for } x \in \Omega_{\rho_{\varepsilon}}.
\end{equation*}
Hence by \eqref{initial behavior for mPorous} of Proposition \ref{prop:initial behavior for mPorous}
there exists $\tau_{2,\varepsilon} \in (0, \tau_{1,\varepsilon}]$ such that for $(x,t) \in \Gamma_{\rho_{\varepsilon}}\times (0, \tau_{2,\varepsilon}]$
\begin{equation}
\label{1st inequality for initial-boundary mPorous}
c(m)\left(1+ {\eta_\varepsilon}\right)^{-\frac 1{1-m}} (\rho_{\varepsilon})^{-\frac 2{1-m}} \le t^{-\frac 1{1-m}}u(x,t) \le 
c(m)\left(1- {\eta_\varepsilon}\right)^{-\frac 1{1-m}} (\rho_{\varepsilon})^{-\frac 2{1-m}},
\end{equation}
where $\Gamma_{\rho_{\varepsilon}}$ is defined by  \eqref{level surfaces of distance functions} with $\rho = \rho_{\varepsilon}$.

Moreover, by \eqref{limits corresponds to boundary blow up sol pm for initial-boundary}, there exists $\tau_{\varepsilon} \in (0, \tau_{2,\varepsilon}]$ such that for $(x,t) \in \Gamma_{\rho_{\varepsilon}}\times (0, \tau_{\varepsilon}]$
\begin{eqnarray*}
&& t^{-\frac 1{1-m}}(\rho_{\varepsilon})^{\frac 2{1-m}} w_{+}(x,t) \ge c(m)\left(1- {\eta_\varepsilon}\right)^{-\frac 1{1-m}},
\\
&&t^{-\frac 1{1-m}}(\rho_{\varepsilon})^{\frac 2{1-m}} w_{-}(x,t) \le c(m)\left(1+ {\eta_\varepsilon}\right)^{-\frac 1{1-m}}.
\end{eqnarray*}
Thus combining these inequalities with \eqref{1st inequality for initial-boundary mPorous} yields that
\begin{equation}
\label{inequalities on Gamma_{epsilon} for ib for mPorous}
w_{-}\le u \le w_{+}\ \mbox{ on } \Gamma_{\rho_{\varepsilon}} \times (0,\tau_{\varepsilon}].
\end{equation}

Observe that
\begin{eqnarray}
&w_{-} < \beta=u <  w_{+}\ &\mbox{ on }\partial\Omega \times (0,\tau_{\varepsilon}],\label{comparison on the boundary original for mPorous}
\\
& w_{-}=u=w_{+}=0\ &\mbox{ on } \Omega_{\rho_{\varepsilon}} \times \{0\}.\label{initial condition comparison for mPorous}
\end{eqnarray}
Therefore, by combining these with \eqref{inequalities on Gamma_{epsilon} for ib for mPorous} and \eqref{differential inequalities for mPorous initial-boundary}, we get the conclusion \eqref{important estimates for the solution near the boundary for short time for mPorous ib} from the comparison principle.
\qed

\vskip 4ex
Let us next consider the solution $u$ of 
 problem {\rm (\ref{cauchy mPorous})}.   Namely, we consider the bounded solution $u=u(x,t)$  of 
  the Cauchy problem:
  \begin{equation*}
\label{cauchy mPorous beta}
\partial_t u= \Delta u^{m}\ \mbox{ in }\ \mathbb R^N \times (0, \infty)\quad\mbox{ and }\ u = \beta\mathcal X_{\Omega^c}\ \mbox{ on }\ \mathbb R^N \times \{0\},
\end{equation*}
where  $\mathcal X_{\Omega^c}$ is the characteristic function of the set $\Omega^c = \mathbb R^N \setminus \Omega$.
Let us set $\phi(s) = s^{m}$ for $s \ge 0$. We use a result from \cite{MSjde2012}:  
for every $\gamma > 0$, there exists a unique $C^2$ solution $f_\gamma = f_\gamma(\xi)$ of the problem:
\begin{eqnarray}
&& \left( \phi^\prime(f_\gamma) f_\gamma^\prime\right)^\prime + \frac 12 \xi  f_\gamma^\prime = 0\ \mbox{ in }\ \mathbb R,
\label{ode selfsimilar cauchy}
\\
&& f_\gamma(-\infty) = \gamma,\quad f_\gamma(\infty) = 0,\label{boundary conditions at infty}
\\
&& f_\gamma^\prime < 0\ \mbox{ in }\ \mathbb R.
\label{monotonicity cauchy}
\end{eqnarray}
Moreover, \cite[Theorem 5 and its example 3, p. 388 and p. 390]{AParma1974} also gives \eqref{limits corresponds to boundary blow up sol for initial-boundary}. Note that, if we put $h(s,t) = f_\gamma\left(t^{-1/2}s\right)$ for $s \in \mathbb R$ and $t > 0$, then $h$ satisfies the one-dimensional problem:
$$
\partial_t h = \partial_s^2 \phi(h)\ \mbox{ in } \mathbb R \times(0,\infty) \mbox{ and } h = \gamma\mathcal X_{(-\infty,0]}\ \mbox{ on } \mathbb R\times\{0\}.
$$
Let $0 < \varepsilon < \frac 14$. By the same proof as in \cite[Proof of (3.35), pp. 251--252]{MSjde2012},
we  find a sufficiently small $0 < \eta_\varepsilon << \varepsilon$ and two $C^2$ functions $f_{\pm} = f_{\pm}(\xi)$ for $\xi  \in \mathbb R$ satisfying:
\begin{eqnarray}
&& f_\pm(\xi) = f_{\beta\pm\varepsilon} \left(\sqrt{1\mp 2\eta_\varepsilon}\ \xi\right)\ \mbox{ if } \xi \ge \eta_\varepsilon, \label{scaling of f on whole line}
\\
&&f_\pm^\prime < 0\ \mbox{ in } \mathbb R,\label{monotonicity both on whole line}
\\
&& f_-(-\infty) < \beta= f_\beta(-\infty) < f_+(-\infty)\ \mbox{ and }\ f_- < f_\beta< f_+\ \mbox{ in } \mathbb R,\label{from above and below on whole line}
\\
&& \left(\phi^\prime(f_\pm)f_\pm^\prime\right)^\prime + \frac 12\xi f_\pm^\prime = h_\pm(\xi) f_\pm^\prime\ \mbox{ in } \mathbb R,\label{ modified ode on whole line}
\end{eqnarray}
and
\begin{equation}
\label{approximation for f_{beta} cauchy}
f_{\pm} \to f_{\beta}\ \mbox{ as } \varepsilon \to 0^{+}\ \mbox{ uniformly on } \mathbb R.
\end{equation}
Moreover, by \eqref{limits corresponds to boundary blow up sol for initial-boundary} we also have
\eqref{limits corresponds to boundary blow up sol pm for initial-boundary}.

By setting
\begin{equation}
\label{super&subsolutions for mPorous c-problem}
w_\pm(x,t) =f_\pm(t^{-1/2}d^{*}(x))\ \mbox{ for } (x,t) \in \mathbb R^{N} \times (0,\infty),
\end{equation}
we obtain
\begin{proposition}
\label{prop:estimates from above and below nbd of boundary by pm for mPorous for cauchy}
Let  $u$ be the solution of 
 problem {\rm (\ref{cauchy mPorous})} where $\partial\Omega$ is bounded and of class $C^2$ and  $g \equiv \beta$ for some positive constant $\beta > 0$. 
 For every small $\varepsilon > 0$ there exist $\rho_\varepsilon \in (0, \rho_0)$ and $\tau_{\varepsilon}>0$ satisfying
\begin{equation}
\label{important estimates for the solution near the boundary for short time for mPorous c}
w_{-} \le u \le w_{+} \ \mbox{ in } \Omega_{\rho_\varepsilon} \times (0, \tau_{\varepsilon}],
\end{equation}
where $w_{\pm}$ are given by \eqref{super&subsolutions for mPorous c-problem} and $\Omega_{\rho_\varepsilon}$ is defined by  \eqref{neighborhood of boundary} with $\rho = \rho_\varepsilon$.
\end{proposition}

\noindent
{\it Proof.}\ The proof is similar to that of Proposition \ref{prop:estimates from above and below nbd of boundary by pm for mPorous}. The ingredient \eqref{comparison on the boundary original for mPorous}  is replaced by   the corresponding inequalities on $\{ x \in \mathbb R^{N} : d^{*}(x) = -\rho_{\varepsilon}\} \times (0,\tau_{\varepsilon}]$.
\qed


\setcounter{equation}{0}
\setcounter{theorem}{0}

\section{Proofs of Theorems \ref{th:interaction curvatures pLaplace} and \ref{th:interaction curvatures mPorous}}
\label{section5}

By virtue of section \ref{section3}, we can assume that $\partial\Omega$ is bounded and of class $C^2$ and  $f \equiv g \equiv \beta$ for some positive constant $\beta > 0$.
We will  use a geometric lemma from \cite{MSprseA2007} adjusted to our situation.

\begin{lemma}{\rm (\cite[Lemma 2.1, p. 376]{MSprseA2007})}
\label{lm:asympvolume}
Let $\kappa_j(y_0)<\frac1{R}$ for every $j=1,\dots,N-1.$ Then we have:
\begin{equation*}
\label{asympvolume}
\lim_{s\to 0^+} s^{-\frac{N-1}{2}} \mathcal H^{N-1}(\Gamma_s\cap B_R(x_0))=2^{\frac{N-1}{2}}\omega_{N-1} \
\left\{\prod_{j=1}^{N-1}\left(\frac1{R}-\kappa_j(y_0)\right)\right\}^{-\frac12},
\end{equation*}
where $\mathcal H^{N-1}$ is the standard $(N-1)$-dimensional Hausdorff measure, and $\omega_{N-1}$ is the volume of the unit ball in $\mathbb R^{N-1}.$
\end{lemma}

Let us first prove Theorem \ref{th:interaction curvatures pLaplace} for the solution $u$ of 
 problem {\rm (\ref{diffusion pLaplace})-(\ref{initial pLaplace})} by using Proposition \ref{prop:estimates from above and below nbd of boundary by pm}.
Take a small $\varepsilon > 0.$  Let $\alpha >  \frac {(N+1)(2-p)}{2p}$. Then Proposition \ref{prop:estimates from above and below nbd of boundary by pm}  yields that for every $t \in (0, \tau_{\varepsilon}]$
\begin{equation}
\label{integral inequalities over the main part of the ball}
\int\limits_{B_R(x_0)\cap\Omega_{\rho_\varepsilon}}\! \left(w_{-}(x,t)\right)^{\alpha}\ dx \le\!\!\!\int\limits_{B_R(x_0)\cap\Omega_{\rho_\varepsilon}}\! \left(u(x,t)\right)^{\alpha}\ dx \le \!\!\!\int\limits_{B_R(x_0)\cap\Omega_{\rho_\varepsilon}}\! \left(w_{+}(x,t)\right)^{\alpha}\ dx.
\end{equation}
For $(x,t) \in \left(\Omega\setminus \Omega_{\rho_\varepsilon}\right) \times (0,\infty)$, by \eqref{upper bound by separable solution for pLaplace} of Proposition \ref{prop:initial behavior for pLaplace}, we have
\begin{equation}
\label{upper bound over the rest part of the ball}
t^{-\frac{N+1}{2p}}\left(u(x,t)\right)^{\alpha} \le t^{-\frac{N+1}{2p} + \frac \alpha{2-p}} (v(x))^\alpha.
\end{equation}
Therefore, since $\overline{B_R(x_0)} \setminus \Omega_{\rho_\varepsilon}$ is a compact set contained in $\Omega$ and $-\frac{N+1}{2p} + \frac \alpha{2-p} > 0$, we see that
\begin{equation}
\label{the rest tends to zero}
t^{-\frac{N+1}{2p}}\!\!\!\int\limits_{B_R(x_0)\setminus\Omega_{\rho_\varepsilon}}\! \left(u(x,t)\right)^{\alpha}\ dx \to 0\ \mbox{ as } t \to 0^+.
\end{equation}
With the aid of the co-area formula, we have
\begin{eqnarray*}
&&\int\limits_{B_R(x_0)\cap\Omega_{\rho_\varepsilon}}\! \left(w_{\pm}(x,t)\right)^{\alpha}\ dx 
\\
&&\qquad\quad= t^{\frac {N+1}{2p}}\int_0^{\rho_\varepsilon t^{-1/2}} \left(\varphi_\pm(\xi)\right)^\alpha
\xi^{\frac{N-1}2}\left(t^{1/p}\xi\right)^{-\frac {N-1}2}\mathcal H^{N-1}\left(B_R(x_0)\cap \Gamma_{t^{1/p}\xi}\right) d\xi.
\end{eqnarray*}
Thus, when $\kappa_j(y_0)<\frac1{R}$ for every $j=1,\dots,N-1,$ 
by Lebesgue's dominated convergence theorem and Lemma \ref{lm:asympvolume},  we get
$$
\lim_{t\to 0^+} t^{-\frac{N+1}{2p}}\!\!\!\int\limits_{B_R(x_0)\cap\Omega_{\rho_\varepsilon}}\!\! (w_\pm)^\alpha\ dx = 2^{\frac{N-1}2}\omega_{N-1}\left\{\prod_{j=1}^{N-1}\left(\frac 1R-\kappa_j(y_0)\right)\right\}^{-\frac 12}\int_0^\infty \left(\varphi_\pm(\xi)\right)^\alpha \xi^{\frac{N-1}2} d\xi.
$$
Here \eqref{limits corresponds to to boundary blow up sol pm} together with the inequality $-\frac {p\alpha}{2-p} + \frac {N-1}2 < -1$ guarantees that the right-hand side of this formula is finite.
Moreover,  by Lebesgue's dominated convergence theorem and \eqref{approximation for varphi}, we see  that
$$
\lim_{\varepsilon \to 0}\int_0^\infty \left(\varphi_\pm(\xi)\right)^\alpha \xi^{\frac{N-1}2} d\xi = \int_0^\infty \left(\varphi(\xi) \right)^\alpha\xi^{\frac{N-1}2} d\xi.
$$
Therefore, since $\varepsilon > 0$ is arbitrarily small,  it follows from \eqref{integral inequalities over the main part of the ball} and \eqref{the rest tends to zero} that \eqref{asymptotics and curvatures pLaplace} holds true, where we set
$$
c = 2^{\frac{N-1}2}\omega_{N-1}\int_0^\infty \left(\varphi(\xi)\right)^\alpha \xi^{\frac{N-1}2} d\xi.
$$
\par
It remains to consider the case where $\kappa_j(y_0) = \frac 1R$ for some $j \in \{ 1, \cdots, N-1\}$. Choose a sequence of balls $\{ B_{R_k}(x_k) \}_{k=1}^\infty$  satisfying:
$$
 R_k < R,\ y_0 \in \partial B_{R_k}(x_k), \mbox{ and } B_{R_k}(x_k) \subset B_R(x_0)\mbox{ for every } k \ge 1,\ \mbox{ and }\ \lim_{k\to \infty}R_k = R.
$$
Since $\kappa_j(y_0)\le \frac 1R<\frac1{R_k}$ for every $j=1,\dots,N-1$ and every $k \ge 1$, we can apply the previous case to each $B_{R_k}(x_k)$ to see
 that for every $k \ge 1$
\begin{eqnarray*}
\liminf_{t \to 0^+}t^{-\frac{N+1}{2p} }\!\!\!\int\limits_{B_R(x_0)} \left(u(x,t)\right)^\alpha \ dx &\ge& \liminf_{t \to 0^+}t^{-\frac{N+1}{2p} }\!\!\!\int\limits_{B_{R_k}(x_k)} \left(u(x,t)\right)^\alpha\ dx 
\\
&=& c\left\{\prod\limits_{j=1}^{N-1}\left(\frac 1{R_k} - \kappa_j(y_0)\right)\right\}^{-\frac 12}.
\end{eqnarray*}
Hence, letting $k \to \infty$ yields that 
$$
\liminf_{t \to 0^+}t^{-\frac{N+1}{2p} }\!\!\!\int\limits_{B_R(x_0)} \left(u(x,t)\right)^\alpha\ dx = \infty,
$$
which completes the proof for problem  {\rm (\ref{diffusion pLaplace})-(\ref{initial pLaplace})}.

\vskip 5ex
The proof of Theorem \ref{th:interaction curvatures pLaplace} for  
 problem {\rm (\ref{cauchy pLaplace})} runs similarly with the aid of Proposition \ref{prop:estimates from above and below nbd of boundary by pm for cauchy}.
 Also, the proof of Theorem \ref{th:interaction curvatures mPorous}  runs similarly with the aid of Propositions \ref{prop:estimates from above and below nbd of boundary by pm for mPorous} and \ref{prop:estimates from above and below nbd of boundary by pm for mPorous for cauchy}. Of course, for problems {\rm (\ref{diffusion mPorous})-(\ref{initial mPorous})} and  {\rm (\ref{cauchy mPorous})},
 we use Proposition \ref{prop:initial behavior for mPorous} and the assumption that $\alpha >  \frac {(N+1)(1-m)}{4}$  instead of Proposition \ref{prop:initial behavior for pLaplace} and the assumption that $\alpha >  \frac {(N+1)(2-p)}{2p}$.


\vskip 4ex
\bigskip
\noindent{\large\bf Acknowledgement.}
\smallskip
The results in the case where $f\equiv g \equiv 1$ were announced in \cite{Soberwolfach2013}. After the author's talk in the Mini-Workshop titled ``The p-Laplacian Operator and Applications'' (10 Feb - 16 Feb 2013) the present paper was completed.
The author gratefully acknowledges the hospitality of the Mathematisches Forschungsinstitut Oberwolfach.

\end{document}